\documentclass[reqno,8pt]{amsart}
\usepackage{amsfonts}

\usepackage{amsmath}

\newtheorem{theorem}{Theorem}
\theoremstyle{plain}

\newtheorem{corollary}{Corollary}

\newtheorem{lemma}{Lemma}

\newtheorem{remark}{Remark}

\numberwithin{equation}{section}

\input{tcilatex}

\begin{document}
\title[$p$-adic Dedekind and Hardy-Berndt Type Sums...]{$p$-adic Dedekind
and Hardy-Berndt Type Sums related to Volkenborn Integral on $\mathbb{Z}_{p}$%
}
\author{YILMAZ\ SIMSEK}
\address{Akdeniz University Faculty of Science Depart. of Mathematics 07058
Antalya-Turkey}
\email{ ysimsek@akdeniz.edu.tr}

\begin{abstract}
The purpose of this paper is to construct $p$-adic Dedekind sums and
Hardy-Berndt type sums. We also construct generating function of the twisted
Bernoulli polynomials and functions. Furthermore, we give some discussions
on elliptic analogue of the Apostol-Dedekind sums.
\end{abstract}

\maketitle

\section{Introduction, Definitions and Notations}

Let $(h,k)=1$ with $k>0$. The classical Dedekind sums are defined by%
\begin{equation*}
s\left( h,k\right) =\sum_{a=1}^{k-1}\left( \left( \frac{a}{k}\right) \right)
\left( \left( \frac{ha}{k}\right) \right) ,
\end{equation*}%
where $\left( \left( x\right) \right) =x-[x]_{G}-\frac{1}{2}$, if $x\notin
\mathbb{Z}$, $\left( \left( x\right) \right) =0$, $x\in \mathbb{Z}$, where $%
[x]_{G}$ is the largest integer $\leq x$ cf. (\cite{1}, \cite{19}, \cite{11}%
, \cite{10}, \cite{simjnt2003dede}, \cite{simjmaaqDed}).

In this paper, $\mathbb{Z}_{p}$, $\mathbb{Q}_{p}$, $\mathbb{C}_{p}$, $%
\mathbb{C}$ and $\mathbb{Z}$, respectively, denote the ring of $p$-adic
integers, the field of $p$-adic rational numbers, the $p$-adic completion of
the algebraic closure of $\mathbb{Q}_{p}$ normalized by $\left| p\right|
_{p}=p^{-1}$, and the complex field and integer numbers. Let $q$ be an
indeterminate such that if $q\in \mathbb{C}$, then $\left| q\right| <1$ and
if $q\in \mathbb{C}_{p}$, then $\left| 1-q\right| _{p}<p^{-1/\left(
p-1\right) }$, so that $q^{x}=$exp$\left( x\text{log}_{p}q\right) $ for $%
\left| x\right| _{p}\leq 1$. Let $\left[ x\right] =\left[ x:q\right] =\frac{%
1-q^{x}}{1-q}$. We note that $\lim_{q\rightarrow 1}\left[ x\right] =x$. The $%
p$-adic $q$-Volkenborn integral is originally constructed by Kim \cite%
{kimJMAA2007}, \cite{MODIFkim}, which is defined as follows: for $g\in UD(%
\mathbb{Z}_{p},\mathbb{C}_{p})=\left\{ g\mid g:\mathbb{Z}_{p}\rightarrow
\mathbb{C}_{p}\text{ is uniformly differentiable function}\right\} $, the $p$%
-adic $q$-Volkenborn integral is defined by $I_{q}\left( g\right) =\int_{%
\mathbb{Z}_{p}}g\left( x\right) d\mu _{q}\left( x\right) =\lim_{N\rightarrow
\infty }\frac{1}{\left[ p^{N}:q\right] }\sum_{x=0}^{p^{N}-1}g\left( x\right)
q^{x}$, where $\mu _{q}\left( x+p^{N}\mathbb{Z}_{p}\right) =\frac{q^{x}}{%
\left[ p^{N}:q\right] }$. Note that $I_{1}\left( g\right)
=\lim_{q\rightarrow 1}I_{q}\left( g\right) $. If $g_{1}\left( x\right)
=g\left( x+1\right) $, then%
\begin{equation}
I_{1}\left( g_{1}\right) =I_{1}\left( g\right) +g^{\prime }\left( 0\right)
\text{ cf. (\cite{kimJMAA2007}, \cite{KimNewApproc}),}  \label{1,1}
\end{equation}%
\begin{equation}
qI_{q}(g_{1})-I_{q}(g)=(q-1)g(0)+\frac{q-1}{\log q}g^{\prime }\left(
0\right) \text{ cf. \cite{MODIFkim},}  \label{ayeni-0}
\end{equation}%
where $g^{\prime }\left( 0\right) =\left. \frac{d}{dx}g\left( x\right)
\right| _{x=0}$. The $q$-deformed $p$-adic invariant integral on $\mathbb{Z}%
_{p}$, in the fermionic sense, is defined by%
\begin{equation}
I_{-1}(f)=\lim_{q\rightarrow -1}I_{q}(f)=\int_{\mathbb{Z}_{p}}f(x)d\mu
_{-1}(x)=\lim_{N\rightarrow \infty }\sum_{x=0}^{p^{N}-1}(-1)^{x}f(x)\text{
cf. (\cite{kimJMAA2007}, \cite{MODIFkim}).}  \label{Equ-18}
\end{equation}

In \cite{kimJMAA2007}, \cite{MODIFkim}, by busing $p$-adic $q$-integral on $%
\mathbb{Z}_{p}$, Kim defined%
\begin{equation}
I_{-1}(f_{1})+I_{-1}(f)=2f(0)\text{.}  \label{Equ-19}
\end{equation}

For applications of the $p$-adic $q$-integral on $\mathbb{Z}_{p}$ see also
cf. (\cite{cenkcisimsekcank}, \cite{30}).

\section{$p$-adic Hardy-Berndt type sums}

In this section, by using the $q$-deformed $p$-adic invariant integral on $%
\mathbb{Z}_{p}$, in the fermionic sense, we construct $p$-adic Hardy-Berndt
type sums.

\begin{equation}
\int_{\mathbb{Z}_{p}}\sin (bx)d\mu _{-1}(x)=-\tan (\frac{b}{2})\text{ cf. %
\cite{kimJMAA2007}.}  \label{ayeni-2}
\end{equation}%
Multiplying both sides of (\ref{ayeni-2}) by $\frac{1}{n\pi }$, and
replacing $b$ by $2\pi ny$, and\ then summing over $n=1,2,...,\infty $, we
have%
\begin{equation*}
\frac{1}{\pi }\sum_{n=1}^{\infty }\frac{1}{n}\int_{\mathbb{Z}_{p}}\sin (2\pi
nyx)d\mu _{-1}(x)=-\frac{1}{\pi }\sum_{n=1}^{\infty }\frac{1}{n}\tan (\pi ny)
\end{equation*}%
After some elementary calculations, we get%
\begin{equation*}
\int_{\mathbb{Z}_{p}}\left( \frac{1}{\pi }\sum_{n=1}^{\infty }\frac{1}{n}%
\sin (2\pi nyx)\right) d\mu _{-1}(x)=-\frac{1}{\pi }\sum_{n=1}^{\infty }%
\frac{1}{n}\tan (\pi ny),
\end{equation*}%
where $\left( \left( y\right) \right) =\frac{1}{\pi }\sum_{n=1}^{\infty }%
\frac{\sin (2\pi ny)}{n}$ cf. (\cite{BG1}, \cite{cenkcisimsekcank}, \cite%
{simjnt2003dede}, \cite{simsekJKMShardy}).Thus we arrive at the following
result:

\begin{lemma}
\label{lem-1}%
\begin{equation*}
\int_{\mathbb{Z}_{p}}\left( \left( yx\right) \right) d\mu _{-1}(x)=\frac{1}{%
\pi }\sum_{n=1}^{\infty }\frac{1}{n}\tan (\pi ny).
\end{equation*}
\end{lemma}

By using Lemma \ref{lem-1}, we construct $p$-adic Hardy-Berndt type sums as
follows:

\begin{theorem}
\label{teo-1}Let $h,k\in \mathbb{Z}$, $(h,k)=1$. If $h$ is odd and $k$ is
even, then we have%
\begin{equation*}
S_{2}(h,k)=-\frac{1}{2}\int_{\mathbb{Z}_{p}}\left( \left( \frac{hx}{k}%
\right) \right) d\mu _{-1}(x).
\end{equation*}
\end{theorem}

\begin{proof}
In\cite{BG1}, Bernd and Goldberg defined Hardy sums, $S_{2}(h,k)$ as
follows: if $h$ is odd and $k$ is even, then
\begin{equation}
S_{2}(h,k)=-\frac{1}{2\pi }\sum_{%
\begin{array}{c}
n=1 \\
2n\not\equiv 0(\func{mod}k)%
\end{array}%
}^{\infty }\frac{\tan (\frac{\pi hn)}{k})}{n}.  \label{Eq-3}
\end{equation}%
By using Lemma \ref{lem-1} with $y=\frac{h}{k}$ with $(h,k)=1$, we have
\begin{equation}
\int_{\mathbb{Z}_{p}}\left( \left( \frac{hx}{k}\right) \right) d\mu _{-1}(x)=%
\frac{1}{\pi }\sum_{n=1}^{\infty }\frac{1}{n}\tan (\frac{\pi hn}{k}).
\label{ayeni-1}
\end{equation}%
By using (\ref{Eq-3}) and (\ref{ayeni-1}), with $2n\not\equiv 0(\func{mod}k)$%
, after some elementary calculations, we arrive at the desired result.
\end{proof}

In\cite{BG1}, Bernd and Goldberg defined Hardy sums, $S_{3}(h,k)$ as
follows: if $k$ is odd, then
\begin{equation}
S_{3}(h,k)=\frac{1}{\pi }\sum_{n=1}^{\infty }\frac{\tan (\frac{\pi hn)}{k})}{%
n}.  \label{Eq-4}
\end{equation}%
By using (\ref{Eq-4}) and Theorem \ref{teo-1}, we obtain the following
corollary:

\begin{corollary}
Let $h,k\in \mathbb{Z}$, $(h,k)=1$. If $k$ is odd, then we have%
\begin{equation*}
S_{3}(h,k)=\int_{\mathbb{Z}_{p}}\left( \left( \frac{hx}{k}\right) \right)
d\mu _{-1}(x).
\end{equation*}
\end{corollary}

Multiplying both sides of (\ref{ayeni-2}) by $\frac{4}{(2n-1)\pi }$, and
replacing $b$ by $\frac{\pi h(2n-1)}{2k}$, with $(h,k)=1$, and\ then summing
over $n=1,2,...,\infty $, we have%
\begin{equation*}
\frac{4}{\pi }\sum_{n=1}^{\infty }\frac{1}{2n-1}\int_{\mathbb{Z}_{p}}\sin (%
\frac{\pi h(2n-1)x}{2k})d\mu _{-1}(x)=-\frac{4}{\pi }\sum_{n=1}^{\infty }%
\frac{1}{2n-1}\tan (\frac{\pi h(2n-1)}{2k})
\end{equation*}%
After some elementary calculations, we get%
\begin{equation}
\int_{\mathbb{Z}_{p}}\left( \frac{4}{\pi }\sum_{n=1}^{\infty }\frac{1}{2n-1}%
\sin (\frac{\pi h(2n-1)x}{2k})\right) d\mu _{-1}(x)=-\frac{4}{\pi }%
\sum_{n=1}^{\infty }\frac{1}{2n-1}\tan (\frac{\pi h(2n-1)}{2k}),
\label{ayeni-3}
\end{equation}%
where%
\begin{equation*}
(-1)^{[x]_{G}}=\frac{4}{\pi }\sum_{n=1}^{\infty }\frac{\sin \left( (2n-1)\pi
x\right) }{2n-1}\text{ cf. (\cite{BG1}, \cite{cenkcisimsekcank}, \cite%
{simjnt2003dede}, \cite{simsekJKMShardy}).}
\end{equation*}%
If $h$ and $k$ are odd, then we have%
\begin{equation}
S_{5}(h,k)=\frac{2}{\pi }\sum_{%
\begin{array}{c}
n=1 \\
2n-1\not\equiv 0(\func{mod}k)%
\end{array}%
}^{\infty }\frac{\tan (\frac{\pi h(2n-1)}{2k})}{2n-1}\text{ cf. \cite{BG1}.}
\label{Eq-6}
\end{equation}%
By using (\ref{ayeni-3}) and (\ref{Eq-6}), we arrive at the following
theorem:

\begin{theorem}
Let $h,k\in \mathbb{Z}$, $(h,k)=1$. If $k$ and $k$ are odd, then we have%
\begin{equation*}
S_{5}(h,k)=-2\int_{\mathbb{Z}_{p}}(-1)^{[\frac{hx}{2k}]_{G}}d\mu _{-1}(x).
\end{equation*}
\end{theorem}

Let $h$ and $k$ denote relatively prime integers with $k>0$. If $h+k$ is
odd, then
\begin{equation}
S(h,k)=\frac{4}{\pi }\sum_{n=1}^{\infty }\frac{\tan (\frac{\pi h(2n-1)}{2k})%
}{2n-1}\text{ cf. \cite{BG1}.}  \label{Eq-1}
\end{equation}

By using (\ref{ayeni-3}) and (\ref{Eq-1}), we easily arrive at the following
corollary:

\begin{corollary}
Let $h,k\in \mathbb{Z}$, $(h,k)=1$ with $k>0$. If $h+k$ is odd, then%
\begin{equation*}
S(h,k)=-\int_{\mathbb{Z}_{p}}(-1)^{[\frac{hx}{2k}]_{G}}d\mu _{-1}(x).
\end{equation*}
\end{corollary}

Note that for detail on Hardy-Berndt sums see also cf. (\cite{BG1}, \cite%
{simjnt2003dede}, \cite{simsekJKMShardy}, \cite{cenkcisimsekcank}).

\section{Twisted Dedekind sums}

In this section, by using $q$-Volkenborn integral, we construct a new
generating function of twisted $q$-Bernoulli polynomials. We define twisted
new approach $q$-Bernoulli functions. We also construct $p$-adic twisted $q$%
-Dedekind type sums. Let $T_{p}=\bigcup\limits_{n\geq 1}C_{p^{n}}=\underset{%
n\rightarrow \infty }{\text{lim}}C_{p^{n}},$where $C_{p^{n}}=\left\{
w:w^{p^{n}}=1\right\} $ is the cyclic group of order $p^{n}$. For $w\in
T_{p} $, the function $x\mapsto w^{x}$ is a locally constant function from $%
\mathbb{Z}_{p}$ to $\mathbb{C}_{p}$ cf.(\cite{kimJMAA2007}, \cite{30}). If
we take $f(x)=w^{x}q^{x}e^{tx}$ in (\ref{ayeni-0}), then we define twisted $%
q $-Bernoulli numbers by means of the following generating function:%
\begin{equation*}
F_{q,w}(t)=I_{q}(w^{x}q^{x}e^{tx})=\frac{(q-1)}{\log q}\frac{(\log q^{2}+t)}{ wq^{2}e^{t}-1}%
=\sum_{n=0}^{\infty }b_{n,w}^{\ast }(q)\frac{t^{n}}{n!},
\end{equation*}%
where the numbers $b_{n,w}^{\ast }(q)$\ are called twisted $q$-Bernoulli
numbers. By using Taylor series of $e^{tz}$ in the above, we get%
\begin{equation}
b_{n,w}^{\ast }(q)=\int_{\mathbb{Z}_{p}}w^{x}q^{x}x^{n}d\mu _{q}(x).
\label{aa0}
\end{equation}%
We define twisted $q$-Bernoulli polynomials by means of the following
generating function:%
\begin{equation}
F_{q,w}(t,z)=F_{q,w}(t)e^{zt}=\sum_{n=0}^{\infty }b_{n,w}^{\ast }(z,q)\frac{%
t^{n}}{n!},  \label{aaa0}
\end{equation}%
where the numbers $b_{n,w}^{\ast }(z,q)$\ are called twisted $q$-Bernoulli
polynomials. By using Cauchy product in the above, we have%
\begin{equation*}
b_{n,w}^{\ast }(z,q)=\int_{\mathbb{Z}_{p}}w^{x}q^{x}(x+z)^{n}d\mu
_{q}(x)=\sum_{k=0}^{n}\left(
\begin{array}{c}
n \\
k%
\end{array}%
\right) z^{n-k}b_{k,w}^{\ast }(q).
\end{equation*}

Observe that if $q\rightarrow 1$, then (\ref{ayeni-0}) reduces to (\ref{1,1}%
). See also cf. (\cite{KimNewApproc}, \cite{30}, \cite{cenkcisimsekcank}).

We need the following definitions. $\overline{B}_{n}\left( x\right) $ is
denoted the $n$th Bernoulli function, which is defined as follows:%
\begin{equation}
\overline{B}_{n}\left( x\right) =\left\{
\begin{array}{ll}
B_{n}(\left\{ x\right\} )\text{,} & \text{if }x\text{ is not an integer} \\
0\text{,} & \text{if }x\text{ is an integer and }n=1\text{,}%
\end{array}%
\right.  \label{A2}
\end{equation}%
where $\left\{ x\right\} $ denotes the fractional part of a real number $x$,
$B_{n}\left( x\right) $ is the Bernoulli polynomial or $\overline{B}%
_{n}\left( x\right) =B_{n}\left( x-\left[ x\right] _{G}\right) $ cf. (\cite%
{1}, \cite{simjmaaqDed}, \cite{simjnt2003dede}, \cite{simsekJKMShardy}).
Apostol \cite{1} generalized Dedekind sums $s\left( h,k,n\right) $ as follows%
\begin{equation}
s(h,k,n)=\sum\limits_{a=1}^{k-1}\frac{a}{k}\overline{B}_{n}\left( \frac{ha}{k%
}\right) ,  \label{A1}
\end{equation}%
where $n,h,k$ are positive integers. By using (\ref{aa0}), (\ref{aaa0}) and (%
\ref{A2}), we have%
\begin{equation}
\overline{b^{\ast }}_{n,w}\left( \frac{jh}{k},q\right) =\int_{\mathbb{Z}%
_{p}}w^{x}q^{x}\left( x+\left\{ \frac{jh}{k}\right\} \right) ^{n}d\mu
_{q}(x).  \label{ayeni-4}
\end{equation}%
By using (\ref{ayeni-4}) and (\ref{A1}), we construct twisted $p$-adic $q$%
-higher order Dedekind type sums by the following theorem:

\begin{theorem}
Let $h,k\in \mathbb{Z}$, $(h,k)=1$, and let $p$ be an odd prime such that $%
p|k$. For $w\in T_{p}$, we have%
\begin{equation*}
s_{w}\left( h,k,m,q\right) =\sum_{j=0}^{k-1}\frac{j}{k}\int_{\mathbb{Z}%
_{p}}w^{x}q^{x}\left( x+\left\{ \frac{jh}{k}\right\} \right) ^{n}d\mu
_{q}(x),
\end{equation*}%
or%
\begin{equation*}
s_{w}\left( h,k,m,q\right) =\sum_{j=0}^{k-1}\frac{j}{k}\sum_{k=0}^{n}\left(
\begin{array}{c}
n \\
k%
\end{array}%
\right) \left\{ \frac{jh}{k}\right\} ^{n-k}b_{k,w}^{\ast
}(q)=\sum_{j=0}^{k-1}\frac{j}{k}\overline{b^{\ast }}_{n,w}\left( \frac{jh}{k}%
,q\right) .
\end{equation*}
\end{theorem}

Observe that when $q\rightarrow 1$ and $w\rightarrow 1$, the sum $%
s_{w}\left( a,b,m,q\right) $ reduces to $k^{m}s\left( h,k,m+1\right) $ in (%
\ref{A1}). Different type $p$-adic Dedekind and Hardy type sums were defined
see for detail (\cite{11}, \cite{10}, \cite{19},\cite{simjmaaqDed}, \cite%
{simsekJKMShardy}, \cite{cenkcisimsekcank}).

\begin{remark}
Recently, elliptic Apostol-Dedekind sums have studied by many authors in the
dfferent areas. Bayad\cite{bayad2}, constructed multiple elliptic Dedekind
sums as an elliptic analogue of Zagier's sums multiple Dedekind sums. Machide%
\cite{T. Machide} defined elliptic analogue of the generalized
Dedekind-Rademacher sums, which involve an elliptic analogue of the
classical Bernoulli functions.

Find elliptic analogue and reciprocity law of the sum $s_{w}\left(
h,k,m,q\right) $.
\end{remark}

\textbf{Acknowledgement.} \textit{This paper was supported by the Scientific
Research Project Administration Akdeniz University.}

\end{document}